\documentclass[11pt,a4paper]{article}
\usepackage{amsfonts}
\usepackage{amsmath}
\usepackage{amssymb}
\usepackage[english]{babel}
\usepackage[T1]{fontenc}
\usepackage{graphicx}
\usepackage{subfigure}
\usepackage{epsfig}
\usepackage{latexsym}
\usepackage{amstext}
\usepackage{amsthm}
\usepackage{graphics}
\usepackage{enumerate}
\usepackage{paralist}
\usepackage{fullpage}
\usepackage{setspace}
\newtheorem{prop}{Proposition}[section]
\newtheorem{rem}[prop]{Remark}

\newtheorem{theo}[prop]{Theorem}
\newtheorem{cor}[prop]{Corollary}

\numberwithin{equation}{section}

\newcommand{\beq}{\begin{eqnarray}}
\newcommand{\beqq}{\begin{eqnarray*}}
\newcommand{\eeq}{\end{eqnarray}}
\newcommand{\eeqq}{\end{eqnarray*}}

\newenvironment{dedication}
  {\center
  \itshape             
  }

\title{ On the windings of complex-valued Ornstein-Uhlenbeck processes
driven by a Brownian motion and by a Stable process }
\author{ {\sc Stavros Vakeroudis} \thanks{University of Cyprus, Department of Mathematics and Statistics, P.O. Box 20537, CY-1678 Nicosia, Cyprus. E-mail: stavros.vakeroudis@gmail.com, Site: https://svakeroudis.wordpress.com/} }
\date{}
\begin{document}

\maketitle
\begin{dedication}
In memoriam, Marc Yor
\end{dedication}
\begin{abstract}
We deal with a complex-valued Ornstein-Uhlenbeck (OU) process with parameter $\lambda\in\mathbb{R}$
starting from a point different from 0 and the way that it winds around the origin.
The starting point of this paper is the skew product representation for an OU process which is associated to the skew product representation
of its driving planar Brownian motion under a new deterministic time scale.
We present the stochastic differential equations (SDEs)
for the radial and for the winding process. Moreover, we obtain the large time (analogue of Spitzer's Theorem for Brownian motion in the complex plane) and the small time asymptotics for the winding and for the radial
process, and we explore the exit time from a cone for a 2-dimensional OU process.
Some Limit Theorems concerning the angle of the cone (when our process winds in a cone) and the parameter $\lambda$ are also presented.
Furthermore, we discuss the decomposition of the winding process of a complex-valued OU process in "small" and "big" windings,
where, for the "big" windings, we use some results already obtained by Bertoin and Werner in \cite{BeW94},
and we show that only the "small" windings contribute in the large time limit.
Finally, we study the windings of a complex-valued OU process driven by a Stable process
and we obtain similar results for its (well-defined) winding and radial process.
\end{abstract}


\tableofcontents

\section{Introduction}
\renewcommand{\thefootnote}{\fnsymbol{footnote}}

Ornstein-Uhlenbeck (OU) processes -initially introduced in \cite{UhO30} as an improvement to Brownian motion (BM)\footnote[2]{When we write: Brownian motion, we always mean real-valued Brownian motion, starting from 0 and planar or complex BM stands for 2-dimensional Brownian motion.} model
in order to describe the movement of a particle-
appear as a natural model (or the limit process of several models)
used in applications of stochastic processes. A reason for that is the character of OU process, that is the
fact that it is positive recurrent, and it has an invariant probability (Gaussian) measure.
This makes its study different (and easier in a way) than that of (planar) complex-valued BM
which is null recurrent.

In particular, the 2-dimensional (complex-valued) OU process and its windings attracted the attention
of many researchers recently, as it turned out to have many applications, namely in the domains of finance and of biology.
For instance, some financial applications can be found e.g. in \cite{LeS98,BNS01,Pat05}, and
for some recent works in a biological context we refer e.g. to the following: rotation of a planar polymer \cite{VYH11},
application in neuroscience \cite{BaG11,DiG12}, etc. Motivated by these applications,
we study here the 2-dimensional OU processes (driven by a BM or by a Stable process) starting from a point distinct from
the origin, and the way that they wind around it.

We start in Section \ref{secpr} by presenting some preliminaries. We recall well-known properties of OU processes with parameter $\lambda>0$ including the key argument of this paper in Proposition \ref{OUwindings}, that is the skew-product representation of complex-valued OU processes starting from a point different from 0. In particular interest is the elementary representation of
the (well-defined) continuous winding process
as the continuous winding process of its driving planar BM, as proven in Vakeroudis \cite{Vak11,Vak11th}.
We note that some other previous discussions concerning OU processes can also be found in Bertoin-Werner \cite{BeW94}.
In that Section, we further give the stochastic differential equations (SDE)
satisfied by the radial and the angular part of our complex-valued OU
and an analogue of Bougerol's identity in terms of OU processes.

Section \ref{secmain} presents the main results concerning the winding number of complex-valued OU processes. In particular,
we study its small and big time asymptotics. We start with stating and proving the small time asymptotics
which is similar to the BM case (Theorem \ref{OUsmalltime}), followed by the analogue for the radial process (Theorem \ref{radiaOUsmallt}). Then, in Theorem \ref{SpitzerOU}, we obtain Spitzer's analogue which essentially says that the (well defined) continuous winding process associated to our complex-valued OU process of parameter $\lambda>0$, starting from a point different from 0, normalized by $t$, converges in law, when $t\rightarrow\infty$, to a Cauchy variable of parameter $\lambda$. Then, we present again the large time asymptotic analogue Theorem for the radial process. Section \ref{secmain} also includes an additional large time asymptotics result for the winding process
and a remark associated to windings in a time interval.

Section \ref{seclambda} deals with some more asymptotics, involving first, the parameter $\lambda$ (big and small $\lambda$ asymptotics)
and second, the asymptotics for the exit time from a cone of complex-valued OU processes for big and small total angle.
In Section \ref{secbigsmall} we discuss the "big" and "small" windings of OU processes, and we compare it
to the BM case (see e.g. \cite{MeY82,PiY84,PiY86,LGY86,PiY89}). In particular, we obtain that the asymptotic behavior
(when $t\rightarrow\infty$) for "big" and "small" windings is quite different for these processes.
We start our study by a result due to Bertoin and Werner \cite{BeW94} (where they use OU processes in order to approach BM)
concerning the "big" windings process for OU processes and we expand it by discussing the contribution of the "small" windings.
More precisely, contrary to the BM case where this decomposition in "big" and "small" windings is fundamental and both processes
affect its winding both in the large time limit and around several points, for OU processes it is essentially only the "small" windings
that are taken into account, a result stated here as Theorem \ref{theosbwindings}.
Loosely speaking, a reason for that is the fact that OU processes are characterized
by a force "pulling" them towards the origin (thus differ from BM), which keeps them in a small neighborhood around it.
Hence, taking into account that OU processes are (positive) recurrent, they are not leaving far away from their origin
consequently it seems that only the "small" windings affect the winding process and not the "big" windings, when $t\rightarrow\infty$.
This Section finishes by a discussion concerning the so-called "very big" windings of a 2-dimensional OU process  (see e.g. \cite{BeW94}).

Finally, Section \ref{OUL} contains a discussion concerning the windings of
complex-valued Ornstein-Uhlenbeck processes driven by a process with jumps (L\'{e}vy process), and in particular by a Stable process (OUSP)
and its small and large time behavior. More precisely, we obtain a stochastic differential equation
satisfied by its well defined winding process, involving the driving Stable process, and the analogue SDE for its associated radial process.
We finish by a discussion concerning a relation between the exit time from a cone for this OU process
with the associated exit time for its driving $\alpha$-stable process ($\alpha\in(0,2]$),
which allows to obtain similar asymptotic results as in Sections \ref{secmain}, \ref{seclambda} and \ref{secbigsmall}.

\section{Reminder on Ornstein-Uhlenbeck processes}\label{secpr}

\subsection{Notations and basic properties}\label{subsecnot}

We start by giving some notations that will be used in what follows.
In addition, we recall some elementary (well-known) properties, concerning
on the one hand Ornstein-Uhlenbeck processes and, on the other hand,
windings of planar Brownian motion, the latter being necessary in order to
study Ornstein-Uhlenbeck windings. Before starting, we note that when we write
$Z$ we will always refer to complex-valued Ornstein-Uhlenbeck process starting
from a point different from 0 (e.g. $z_{0}\in \mathbb{C}^{\ast}$), whereas $B$
will refer to planar Brownian motion (starting from the same point $z_{0}$).

\subsubsection*{Preliminaries on Ornstein-Uhlenbeck processes}
We consider a complex-valued Ornstein-Uhlenbeck (OU) process
\beq\label{OUequation}
    Z_{t} = z_{0} + W_{t} - \lambda \int^{t}_{0} Z_{s} ds,
\eeq
with $\left(W_{t},t\geq0\right)$ denoting a planar Brownian motion with $W_{0}=0$, $z_{0}\in \mathbb{C}^{\ast}$ and $\lambda \geq 0$.
For OU processes, we consider $\left(B_{t},t\geq0\right)$ another planar Brownian motion starting from $z_{0}$,
and we have the following representation (see e.g. \cite{ReY99})
\beqq
    Z_{t} &=& e^{-\lambda t} \left( z_{0} + \int^{t}_{0} e^{\lambda s} dW_{s} \right) \label{OUeB} \\
          &=& e^{-\lambda t} \left( B_{\alpha_{t}} \right), \label{OUeB1}
\eeqq
where
\beq\label{alpha}
    \alpha_{t}=\int^{t}_{0} e^{2\lambda s} ds = \frac{e^{2\lambda t}-1}{2 \lambda} \ \ ; \ \
    \alpha^{-1}_{s}=\frac{1}{2\lambda} \log\left(1+2\lambda s\right) .
\eeq
Note that the first equation can be easily verified by simply applying
Itô's formula on the right hand side of (\ref{OUeB}) in order to obtain (\ref{OUequation}), and the second one follows by invoking
Dambis-Dubins-Schwarz Theorem which states that there exists a planar BM $B$ such that (\ref{OUeB1}) is satisfied. \\
From now on, for simplicity and without any loss of generality, we may consider: $z_{0}=1+i0$, which is
really no restriction.

\begin{prop}\label{propscaling}
Ornstein-Uhlenbeck processes satisfy the following "scaling type" property: for every $t>0$ fixed and $a>0$,
\beqq
Z_{at}\stackrel{(law)}{=}e^{-\lambda (1+a) t} \sqrt{\frac{e^{2\lambda a t}-1}{e^{2\lambda  t}-1}} \ Z'_{t} ,
\eeqq
where $Z'$ is an independent copy of $Z$.
\end{prop}
\begin{proof}
Starting from (\ref{OUeB1}) and using the scaling property of BM, we have: for $a>0$,
\beqq
Z_{at}=e^{-\lambda at}B_{\alpha_{at}}\stackrel{(law)}{=}e^{-\lambda (1+a)t}\sqrt{\frac{\alpha_{at}}{\alpha_{t}}} \ e^{\lambda t} B'_{\alpha_{t}} ,
\eeqq
with $B'$ denoting an independent copy of $B$. \\
The proof finishes by remarking that $Z'_{t}=e^{\lambda t}B'_{\alpha_{t}}$ and
\beqq
\frac{\alpha_{at}}{\alpha_{t}}=\frac{e^{2\lambda a t}-1}{e^{2\lambda  t}-1} .
\eeqq
\end{proof}

\subsubsection*{Skew-product representation of planar Brownian motion}
Before proceeding to the study of complex-valued OU processes, we first recall some useful
results concerning planar BM $B$ starting from $1+i0$, that we will also use later on.
As $B$ starts from a point different from 0,
the continuous winding process of the planar BM  $B$, namely
$$\theta^{B}_{t}=\mathrm{Im}\left(\int^{t}_{0}\frac{dB_{s}}{B_{s}}\right),t\geq0$$ is well defined \cite{ItMK65}.
We also define the radial process of the planar BM $B$:
$$R^{B}_{t}=|B_{t}|\Longrightarrow \log R^{B}_{t}=\mathrm{Re}\left(\int^{t}_{0}\frac{dB_{s}}{B_{s}}\right), t\geq0.$$
Hence, we recall the well-known skew product representation of planar BM $B$ (see also e.g. \cite{ReY99})
\beq\label{skew-product}
\log\left|B_{t}\right|+i\theta_{t}\equiv\int^{t}_{0}\frac{dB_{s}}{B_{s}}=\left(
\beta_{u}+i\gamma_{u}\right)
\Bigm|_{u=H_{t}=\int^{t}_{0}\frac{ds}{\left|B_{s}\right|^{2}}} ,
\eeq
with $(\beta_{u}+i\gamma_{u},u\geq0)$ denoting another planar
Brownian motion starting from $\log 1+i0=0$. \\
Equivalently, (\ref{skew-product}) can also be stated as
\beq\label{skew-product2}
\log\left|B_{t}\right|=\beta_{H_{t}} \ ; \ \ \theta^{B}_{t}=\gamma_{H_{t}} ,
\eeq
and we easily deduce that the two $\sigma$-fields $\sigma \{\left|B_{t}\right|,t\geq0\}$ and $\sigma \{\beta_{u},u\geq0\}$ are identical, whereas $(\gamma_{u},u\geq0)$ is independent from $(\left|B_{t}\right|,t\geq0)$. Note that the inverse of $H$ will play an essential role in the sequel and is given by (for further study of the Bessel clock $H$, see also \cite{Yor80}):
\beqq
A_{u}\equiv H_{u}^{-1}=\inf\{ t:H_{t}>u \}=\int^{u}_{0}e^{2\beta_{s}}ds.
\eeqq

\subsubsection*{Skew-product representation of Ornstein-Uhlenbeck processes}
We return now to the complex-valued OU process $Z$. Similarly to planar BM,
as $Z$ starts from a point different from the origin, the continuous winding process associated to $Z$:
$$\theta^{Z}_{t}=\mathrm{Im}\left(\int^{t}_{0}\frac{dZ_{s}}{Z_{s}}\right),t\geq0$$
is well defined, and we also introduce the associated radial process:
$$R^{Z}_{t}=|Z_{t}|\Longrightarrow \log R^{Z}_{t}=\mathrm{Re}\left(\int^{t}_{0}\frac{dZ_{s}}{Z_{s}}\right), t\geq0.$$
\begin{prop}\label{OUwindings}
For a complex-valued OU process $Z$ we have the following skew-product representation:
\beq
\theta^{Z}_{t} &=& \gamma_{H_{\alpha(t)}}, \label{thetaBMOUSP} \\
\log R^{Z}_{t} &=& \beta_{H_{\alpha(t)}}-\lambda t, \label{RBMOUSP}
\eeq
where $\alpha_{t}=\frac{e^{2\lambda t}-1}{2 \lambda}$.
\end{prop}
\begin{proof}
It follows directly from (\ref{OUeB1}) together with the skew-product representation of BM \eqref{skew-product2}.
Indeed, recalling from Vakeroudis \cite{Vak11} that \eqref{OUeB1} yields
\beq
\theta^{Z}_{t} &=& \theta^{B}_{\alpha_{t}}, \label{thetaBMOU}
\eeq
we get \eqref{thetaBMOUSP}.
Concerning the radial part, using (\ref{OUeB1}) we heve
\beq
\log R^{Z}_{t} &=& \log R^{B}_{\alpha_{t}}-\lambda t, \label{RBMOU}
\eeq
hence \eqref{RBMOUSP}.
\end{proof}
{\noindent We define now the first exit time from a cone with a single boundary $c>0$ for $B$
(respectively for $Z$)\footnote[3]{Note that in what follows, the index $(\lambda)$ of the hitting times (wherever there is one) will always refer to the respective hitting time of an OU process with parameter $\lambda$.} }
\beq\label{TthetaOU}
T^{\theta}_{c} \equiv \inf\left\{t\geq 0 : \theta^{B}_{t}=c \right\}
\ \ (\mathrm{respectively} \ T^{\theta(\lambda)}_{c} \equiv \inf \left\{t\geq 0 : \theta^{Z}_{t}=c \right\}) .
\eeq
We also define the first exit time from a cone with two symmetric boundaries of equal angles $c>0$ for $B$
(respectively for $Z$)
\beqq
T^{|\theta|}_{c} \equiv \inf\left\{t\geq 0 : \left|\theta^{B}_{t}\right|=c \right\} \ \
(\mathrm{respectively} \ T^{|\theta(\lambda)|}_{c} \equiv \inf\left\{t\geq 0 : \left|\theta^{Z}_{t}\right|=c \right\}).
\eeqq
We remark here that we could also study the first exit time from a cone with two different angles $c>0$
and $d>0$, but, for simplicity, we consider only $c=d$.
\begin{cor}\label{OUT} Using the previously introduced notations, we have
\beq
    T^{\theta(\lambda)}_{c}&=&\frac{1}{2\lambda}\log \left(1+2\lambda T^{\theta}_{c}\right);  \label{Ttheta} \\
    T^{|\theta(\lambda)|}_{c}&=&\frac{1}{2\lambda}\log \left(1+2\lambda T^{|\theta|}_{c}\right).  \label{Tabstheta}
\eeq
\end{cor}
\begin{proof} We prove (\ref{Ttheta}) ((\ref{Tabstheta}) follows by repeating the same arguments for $T^{|\theta(\lambda)|}_{c}$). From (\ref{TthetaOU}) and using (\ref{thetaBMOU}), we have
\beqq
    T^{\theta(\lambda)}_{c} = \inf \left\{t\geq 0 : \theta^{B}_{\alpha_{t}}=c \right\}.
\eeqq
Hence
\beq\label{Tlambda}
    T^{\theta(\lambda)}_{c}=\alpha^{-1}\left(T^{\theta}_{c}\right),
\eeq
with $\alpha^{-1}(t)=\frac{1}{2\lambda}\log \left(1+2\lambda t\right)$, which yields (\ref{Ttheta}).
\end{proof}
\begin{rem}
For several asymptotic results of these exit times from a cone, involving small and large values of the parameter $\lambda$
and the angle $c$, we refer to Section \ref{seclambda} below.
\end{rem}

\subsection{Stochastic differential equations satisfied by the radial and angular part}\label{subsecSDE}
In this Subsection, we investigate the stochastic differential equations (SDE) satisfied by
the radial and the angular parts of complex-valued OU processes. For this, we present two SDEs
for both the radial and the angular process, the first one involving the new time scale $\alpha_{t}$ and the second one
based on the initial SDE (\ref{OUequation}) satisfied by our 2-dimensional OU process.

\subsubsection*{First SDE:}
On the one hand, we remark that (\ref{thetaBMOU}) yields
that the winding process for complex-valued OU processes satisfies the same stochastic differential equation
with that of the winding process for planar BM but with a different diffusion coefficient, depending on $\lambda$.
Indeed, we may write the standard planar Brownian motion as
$\left(B_{t}=B^{(1)}_{t}+iB^{(2)}_{t},t\geq0\right)$ starting from $1+i0$, where $(B^{(1)}_{t},t\geq0)$ and $(B^{(2)}_{t},t\geq0)$
are two independent linear BMs starting respectively from $1$ and $0$. Hence
(following e.g \cite{LGY86} or \cite{ReY99} Theorem 2.11 in Chapter V, p. 193)
\beq\label{SDEradial}
\log |Z_{t}| = \log |B_{\alpha_{t}}|= - \lambda t + \mathrm{Re}\left(\int^{\alpha_{t}}_{0}\frac{dB_{s}}{B_{s}}\right)
= - \lambda t  + \int^{\alpha_{t}}_{0}\frac{B^{(1)}_{s}dB^{(1)}_{s}+B^{(2)}_{s}dB^{(2)}_{s}}{|B_{s}|^{2}} .
\eeq
Similarly
\beq\label{SDE}
\theta^{Z}_{t} = \theta^{B}_{\alpha_{t}}&=& \mathrm{Im}\left(\int^{\alpha_{t}}_{0}\frac{dB_{s}}{B_{s}}\right)
= \int^{\alpha_{t}}_{0}\frac{-B^{(2)}_{s}dB^{(1)}_{s}+B^{(1)}_{s}dB^{(2)}_{s}}{|B_{s}|^{2}} .
\eeq
Equivalently, we have in differential form
\beq\label{SDE2}
d(\log |Z_{t}|) &=& - \lambda \ dt+ \left(\frac{B^{(1)}_{u}}{|B_{u}|^{2}} \ dB^{(1)}_{u}
+\frac{B^{(2)}_{u}}{|B_{u}|^{2}} \ dB^{(2)}_{u}\right)\Bigg|_{u=\alpha_{t}=\frac{e^{2\lambda t}-1}{2\lambda}} ; \\
d\theta^{Z}_{t} &=& \left(\frac{-B^{(2)}_{u}}{|B_{u}|^{2}} \ dB^{(1)}_{u}
+\frac{B^{(1)}_{u}}{|B_{u}|^{2}} \ dB^{(2)}_{u}\right)\Bigg|_{u=\alpha_{t}=\frac{e^{2\lambda t}-1}{2\lambda}} .
\eeq
We also remark that skew product representation (\ref{skew-product2}) follows from (\ref{SDEradial}) and (\ref{SDE})
by Dambis-Dubins-Schwarz Theorem.

\subsubsection*{Second SDE:}
Following \cite{LGY86}, we decompose the processes in (\ref{OUequation}) into their real and imaginary coordinates,
that is: $Z_{t}=Z^{(1)}_{t}+iZ^{(2)}_{t}$ and $W_{t}=W^{(1)}_{t}+iW^{(2)}_{t}$,
where $Z^{(1)}$ and $Z^{(2)}$ are two real-valued OU processes, starting respectively from 1 and 0,
$W^{(1)}$ and $W^{(2)}$ are two real-valued BMs starting both from 0, and all of them are independent.
Hence
\beqq
Z_{t}=Z^{(1)}_{t}+iZ^{(2)}_{t}=|Z_{t}|\exp\left(i \theta^{Z}_{t}\right),
\eeqq
and taking logarithms, we get
\beqq
\log|Z_{t}| +i\theta^{Z}_{t}&=&\log Z_{t}=\int^{t}_{0} \frac{dZ_{s}}{Z_{s}}= \int^{t}_{0} \frac{dW_{s}-\lambda Z_{s} ds}{Z_{s}} \\
&=&  \int^{t}_{0} \frac{dW^{(1)}_{s}+i \ dW^{(2)}_{s}}{Z_{s}}-\lambda t
=  \int^{t}_{0} \frac{dW^{(1)}_{s}+i \ dW^{(2)}_{s}}{Z^{(1)}_{t}+iZ^{(2)}_{t}}-\lambda t ,
\eeqq
thus
\beq
\log|Z_{t}|&=& \int^{t}_{0} \frac{Z^{(1)}_{s} dW^{(1)}_{s}+ Z^{(2)}_{s} dW^{(2)}_{s}}{|Z_{s}|^{2}}-\lambda t ; \nonumber \\
\theta^{Z}_{t}&=& \int^{t}_{0} \frac{-Z^{(2)}_{s} dW^{(1)}_{s}+ Z^{(1)}_{s} dW^{(2)}_{s}}{|Z_{s}|^{2}} , \label{SDEthetaZ}
\eeq
and equivalently, in differential form
\beq
d\left(\log|Z_{t}|\right)&=&  \frac{Z^{(1)}_{t}}{|Z_{t}|^{2}} \ dW^{(1)}_{t}+ \frac{Z^{(2)}_{t}}{|Z_{t}|^{2}} \ dW^{(2)}_{t}-\lambda dt ; \nonumber \\
d\theta^{Z}_{t}&=& \frac{-Z^{(2)}_{t}}{|Z_{t}|^{2}} \ dW^{(1)}_{t}+ \frac{Z^{(1)}_{t}}{|Z_{t}|^{2}} \ dW^{(2)}_{t} . \label{SDEdiftheta}
\eeq
With $<\cdot>$ standing for the quadratic variation, we have
\beqq
<Z^{(1)}>_{t}=<Z^{(2)}>_{t}=<W^{(1)}>_{t}=<W^{(2)}>_{t}=t .
\eeqq
Consider $\left(\delta_{t},t\geq0\right)$, $\left(\hat{\delta}_{t},t\geq0\right)$, $\left(b_{t},t\geq0\right)$ and $\left(\hat{b}_{t},t\geq0\right)$
four real BMs all starting from 0, and independent from each other and from all the other processes.
Hence, invoking Dambis-Dubins-Schwarz Theorem, (\ref{SDEthetaZ}) (or equivalently (\ref{SDEdiftheta})) can also be stated in the following form:
\beqq
\log|Z_{t}|&=&  \hat{\delta}_{\int^{t}_{0} \frac{ds}{|Z_{s}|^{2}}}-\lambda t= \int^{t}_{0} \frac{d\hat{b}_{s}}{|Z_{s}|}-\lambda t ; \\
\theta^{Z}_{t}&=& \delta_{\int^{t}_{0} \frac{ds}{|Z_{s}|^{2}}} = \int^{t}_{0} \frac{db_{s}}{|Z_{s}|} .
\eeqq
Note that the latter is the OU analogue of the one for BM, that is (see e.g. \cite[Chapter IV, equation 35.14]{RoW87})
with an independent real BM, starting from 0,
\beqq
d\theta^{B}_{t}=\frac{1}{|B_{t}|} \ db_{t} .
\eeqq
For a similar discussion, see also \cite[Section 4.4.5]{Gar83}).
\begin{rem}
We remark that the two SDEs (\ref{SDE2}) and (\ref{SDEdiftheta}) associated to the winding process of $Z$
are equivalent. This is clear if we replace each OU process in (\ref{SDEdiftheta})
by its equivalent form involving a BM multiplied by $e^{-\lambda t}$ (like in (\ref{OUeB1})).
\end{rem}

\subsection{An expression related to Bougerol's identity in law}\label{subsecBoug}
We can now present the following Proposition coming from \cite{Vak11} which is essentially an attempt to obtain
an analogue of Bougerol's identity in law for Ornstein-Uhlenbeck processes.
We first recall that Bougerol's celebrated identity in law states that:
with $(\beta_{t},t\geq0)$ and $(\hat{\beta_{t}},t\geq0)$ two real independent BMs, for every $u>0$ fixed,
\beq\label{boug}
\sinh(\beta_{u}) \stackrel{(law)}{=} \hat{\beta}_{A_{u}=(\int^{u}_{0}ds\exp(2\beta_{s}))}.
\eeq
For further details and other equivalent expressions and extensions of (\ref{boug}), we refer the interested reader
to \cite{Vak12} and the references therein.
\begin{prop}\label{OUbougerol}
We consider two independent OU processes: $(Z_{t},t\geq 0)$ which is complex-valued and $(\Xi_{t},t\geq
0)$ which is real-valued OU, both starting from a point different from 0. For every $r>0$, define
$T^{(\lambda)}_{r}(\Xi)= \inf \left\{t\geq 0 : e^{\lambda
t } \Xi_{t}=r \right\}$. Then,
\beq\label{OUbougerolequation}
    \theta^{Z}_{T^{(\lambda)}_{r}(\Xi)}  \stackrel{(law)}{=} C_{a(r)},
\eeq
where $a(x)= \arg \sinh (x)$, and $C_{\sigma}$ is a Cauchy variable with parameter $\sigma$.
\end{prop}
\begin{proof}
First, for a real BM $\beta$, we introduce the hitting time of a level $k>0$: $T^{\beta}_{k}= \inf \left\{t\geq 0 : \beta_{t}=k \right\}$.
Taking equation (\ref{OUeB}) or (\ref{OUeB1}) for $\Xi^{\lambda}_{t}$, we have
\beqq
    e^{\lambda t } \Xi_{t}= \delta_{( \frac{e^{2\lambda t}-1}{2 \lambda} )},
\eeqq
with $(\delta_{t},t\geq 0)$ denoting a real Brownian motion starting from the same point with $\Xi$,
different from 0 (without loss of generality, starting e.g. from 1). Thus:
\beq\label{TbU}
    T^{(\lambda)}_{r}(\Xi)=\frac{1}{2\lambda}\log \left(1+2\lambda T^{\delta}_{r}\right).
\eeq
Equation (\ref{thetaBMOU}) for $t=\frac{1}{2\lambda} \log
\left(1+2\lambda T^{\delta}_{r}\right)$, equivalently
$\alpha_{t}=T^{\delta}_{r}$, becomes
\beqq
\theta^{Z}_{T^{(\lambda)}_{r}(\Xi) } =
\theta^{Z}_{\frac{1}{2\lambda}\log \left(1+2\lambda
T^{\delta}_{r}\right) }= \theta^{B}_{u=T^{\delta}_{r}}.
\eeqq
Invoking the skew-product representation (\ref{skew-product2}), we get
\beqq
\theta^{B}_{T^{\delta}_{r}}=\gamma_{H_{T^{\delta}_{r}}}.
\eeqq
The symmetry principle (see \cite{And87} for the original Note and \cite{Gal08} for a detailed discussion),
yields that Bougerol's identity may be equivalently stated as (the bar stands for the supremum)
\beqq
    \sinh(\bar{\beta}_{u}) \stackrel{(law)}{=} \bar{\delta}_{A_{u}},
\eeqq
hence, by identifying the laws of the first hitting times of a level $r>0$, we obtain: $T^{\beta}_{a(r)} \stackrel{(law)}{=} H_{T^{\delta}_{r}}$.
We point out that $H$ is the inverse of $A$ (see e.g. \cite{Vak11}).
The proof finishes by recalling that
$(\gamma_{T^{\beta}_{u}},u\geq0)$ is equal in law to a Cauchy process $(C_{u},u\geq0)$ .
\end{proof}
\begin{rem}
Equation (\ref{TbU}), yields a simple computation of the Laplace transform of $T^{(\lambda)}_{r}(\Xi)$.
More precisely, for $r>1$ (note that we have supposed that $\Xi_{0}=1$),
\beq\label{TLOU}
E\left[\exp\left(-\mu T^{(\lambda)}_{r}(\Xi)\right)\right]=
\frac{1}{\Gamma\left(\frac{\mu}{2\lambda}\right)} \int^{\infty}_{0} dt \ t^{\frac{\mu}{2\lambda}-1} e^{-t-r\sqrt{2\lambda t}} .
\eeq
Indeed, from (\ref{TbU}), using that $E[\exp(-\mu T^{\delta}_{r})]=\exp(-r \sqrt{2\mu})$ (see e.g. \cite{ReY99}), we have that, for every $\mu>0$,
\beqq
E\left[\exp\left(-\mu T^{(\lambda)}_{r}(\Xi)\right)\right]&=&
E\left[\exp\left(- \frac{\mu}{2\lambda}\log \left(1+2\lambda T^{\delta}_{r}\right)\right)\right] \\
&=& E\left[\left(1+2\lambda T^{\delta}_{r}\right)^{-\mu/(2\lambda)}\right] \\
&=& \frac{1}{\Gamma\left(\frac{\mu}{2\lambda}\right)} \int^{\infty}_{0} dt \ t^{\frac{\mu}{2\lambda}-1}
E\left[\exp\left(-t(1+2\lambda T^{\delta}_{r})\right)\right] ,
\eeqq
from which follows (\ref{TLOU}). \\ \\
We note that a similar formula for the Laplace transform of the first hitting time
$$\hat{T}^{(\lambda)}_{r}(\Xi)= \inf \left\{t\geq 0 : \Xi_{t}=r \right\}$$
can be found e.g. in \cite{BoS02} (Chapter 7, Formula 2.0.1, p. 542) or \cite{APP05} (Proposition 2.1 therein; see also \cite{BeH51,Bre67,Sie51}).
In particular, for $r>1$ (recall that $\Xi_{0}=1$),
\beqq
E\left[\exp\left(-\mu \hat{T}^{(\lambda)}_{r}(\Xi)\right)\right]=\frac{H_{-\mu/\lambda}(-\sqrt{\lambda})}{H_{-\mu/\lambda}(-r\sqrt{\lambda})}
=\frac{e^{\lambda /2}D_{-\mu/\lambda}(-\sqrt{2\lambda})}{e^{\lambda r^{2}/2} D_{-\mu/\lambda}(-r\sqrt{2\lambda})} ,
\eeqq
where $H_{\nu}(\cdot)$ is the Hermite function and $D_{\nu}(\cdot)$ is the parabolic cylinder function.
\end{rem}
\begin{rem}
Taking $\lambda=0$ in (\ref{OUbougerolequation}), we obtain
\beqq
\theta_{T^{\delta}_{r}} \stackrel{(law)}{=} C_{a(r)},
\eeqq
where $T^{\delta}_{r}=\inf\{ t:\delta_{t}=r\}$, which is the corresponding result for
planar BM and which is equivalent to Bougerol's identity (\ref{boug}). For more details, see e.g. \cite{Vak11,Vak12}.
\end{rem}

\section{Small and Large time asymptotics}\label{secmain}

\subsection{Small time asymptotics}
Let us now study the windings of complex-valued OU processes in the small time limit.
Starting from Proposition \ref{OUwindings}, we obtain the following:
\begin{theo}\label{OUsmalltime}
The family of processes $$\left(t^{-1/2} \theta^{Z}_{st},s\geq0\right)$$ converges in distribution, as $t\rightarrow0$, to a 1-dimensional Brownian motion $\left(\gamma_{s},s\geq0\right)$.
\end{theo}
\begin{proof}
We follow the main steps of Theorem 7 in Doney-Vakeroudis \cite[p. 297]{DoV12} and we also make use of (\ref{thetaBMOU}).
We split the proof into two parts. \\ \\
$\left.\mathrm{i}\right)$ First, we prove that for the clock $H_{t}=\int^{t}_{0} |B_{s}|^{-2}ds$,
associated to the planar BM $B$ from (\ref{skew-product}) or (\ref{skew-product2}), we have the a.s. convergence:
\beq\label{clocktsmall}
\left(\frac{H(x\alpha_{u})}{\alpha_{u}}, x\geq 0\right) \overset{{a.s.}}{\underset{u\rightarrow0}\longrightarrow} \left(x, x\geq0\right) ,
\eeq
which also implies the weak convergence in the sense of Skorokhod ("$\Longrightarrow$" denotes this type of convergence):
\beq\label{clockBM2}
\left(\frac{H(x\alpha_{u})}{\alpha_{u}}, x\geq 0\right) \overset{{(d)}}{\underset{u\rightarrow0}\Longrightarrow} \left(x, x\geq0\right) .
\eeq
Indeed, using the definition of $H$, we have
\beqq
\frac{H(x\alpha_{u})}{\alpha_{u}}=\frac{1}{\alpha_{u}} \int^{x\alpha_{u}}_{0} \frac{ds}{|B_{s}|^{2}} .
\eeqq
Hence, for every $x_{0}>0$, because $|B_{u}|^{2}\overset{{a.s.}}{\underset{u\rightarrow0}\longrightarrow} 1$,
\beq\label{prsmall1}
\sup_{x\leq x_{0}} \left|\frac{H(x\alpha_{u})-x\alpha_{u}}{\alpha_{u}}\right|&=&
\sup_{x\leq x_{0}} \frac{1}{\alpha_{u}}\left|\int^{x\alpha_{u}}_{0} \left(\frac{1}{|B_{s}|^{2}}-1\right) ds \ \right|
\leq \frac{1}{\alpha_{u}}\int^{x_{0}\alpha_{u}}_{0} \left|\frac{1}{|B_{s}|^{2}}-1\right| ds \nonumber \\
&\stackrel{s=\alpha_{u}w}{=}& \int^{x_{0}}_{0} \left|\frac{1}{|B_{\alpha_{u}w}|^{2}}-1\right| dw \
\overset{{a.s.}}{\underset{u\rightarrow0}\longrightarrow} 0 .
\eeq
Hence, as (\ref{prsmall1}) is true for every $x_{0}>0$, we obtain (\ref{clocktsmall}), thus also (\ref{clockBM2}). \\
Note that this argument is also valid for a more general clock than that of BM. We just have to replace the order of stability
(power 2 in the denominator) by the new order of stability in $\left(\right.0,2\left.\right]$ (for further details see \cite{DoV12}).
\\ \\
$\left.\mathrm{ii}\right)$ Using the skew product representation (\ref{skew-product2})
and the scaling property of BM, we have that for every $s>0$,
\beqq
t^{-1/2} \theta^{Z}_{st}= t^{-1/2} \theta^{B}_{\alpha_{st}}&=& t^{-1/2}\gamma_{\left(H_{\alpha_{st}}\right)}\stackrel{(law)}=\gamma_{\left(t^{-1}H_{\alpha_{st}}\right)} \nonumber \\
&=& \gamma_{\left(\frac{\alpha(st)}{t}\frac{H_{\alpha(st)}}{\alpha(st)}\right)} .
\eeqq
However, we have that
\beq\label{alphaproof}
\frac{\alpha(st)}{t}=\frac{e^{2\lambda st}-1}{2\lambda t}\stackrel{t\rightarrow0}{\longrightarrow} s ,
\eeq
which, together with (\ref{clocktsmall}), finishes the proof.
\end{proof}
{\noindent For the small time limit of the radial process of an Ornstein-Uhlenbeck process, that is: $R^{Z}=(R^{Z}_{u}, u\geq0)=(|Z_{u}|, u\geq0)$, we have:}
\begin{theo}\label{radiaOUsmallt}
The family of processes $$\left(t^{-1/2} \log R^{Z}_{st},s\geq0\right)$$ converges in distribution, as $t\rightarrow0$, to a 1-dimensional Brownian motion $\left(\beta_{s},s\geq0\right)$.
\end{theo}
\begin{proof}
Our proof follows the lines of the proof of Theorem \ref{OUsmalltime}.
Using \eqref{RBMOUSP}, we get:
$$t^{-1/2}\log R^{Z}_{st} = t^{-1/2}\beta_{H_{\alpha(st)}}-\lambda st^{1/2}.$$
The scaling property of BM yields that for every $s>0$,
$$t^{-1/2}\beta_{H_{\alpha(st)}}\stackrel{(law)}{=} \beta_{\left(\frac{\alpha(st)}{t} \frac{H_{\alpha(st)}}{\alpha(st)}\right)}.$$
The proof finishes by invoking \eqref{alphaproof} and the a.s. convergence (\ref{clocktsmall}) of the clock $H$.
\end{proof}

\subsection{Large time asymptotics}\label{largetimeasymp}
Now we turn our study to the Large time asymptotics of the winding process associated to complex-valued OU processes.
Before starting, let us first recall the well-known Spitzer's celebrated asymptotic Theorem for planar BM \cite{Spi58} stating that
\beq\label{Spi}
\frac{2}{\log t} \; \theta^{B}_{t} \overset{{(law)}}{\underset{t\rightarrow\infty}\longrightarrow} C_{1} .
\eeq
For other proofs of this Theorem, the interested reader is refereed to e.g. \cite{Wil74,Dur82,MeY82,BeW94,Yor97,Vak11,VaY12} etc.
Note also that, in a more general framework, the asymptotic behavior of the well defined winding process $\vartheta$
of a planar diffusion starting from a point different from the origin has been discussed by Friedman-Pinsky in \cite{FrP73,FrP73a}
and they showed that, when $t\rightarrow\infty$, $\vartheta_{t}/t$ exists a.s. under some assumptions
meaning that the process winds asymptotically around a point.
For other similar studies, see also Le Gall-Yor \cite{LGY86}.
\\ \\
The following is the analogue of Spitzer's Theorem for OU processes:
\begin{theo}\label{SpitzerOU} \emph{\textbf{(Spitzer's Theorem for OU processes) }} \\
The following convergence in law holds:
\beq\label{SpitzerOUcvg}
\frac{\theta^{Z}_{t}}{t} \;  \overset{{(law)}}{\underset{t\rightarrow\infty}\longrightarrow} C_{\lambda} ,
\eeq
where we recall that, $C_{\sigma}$ is a Cauchy variable with parameter $\sigma$.
\end{theo}
\begin{proof} Using (\ref{thetaBMOU}), we have
\beqq
\frac{\theta^{Z}_{t}}{\lambda t}=\frac{\theta^{B}_{\alpha_{t}}}{\lambda t}
=\frac{\log \alpha_{t}}{2\lambda t} \ \frac{2\theta^{B}_{\alpha_{t}}}{\log \alpha_{t}} .
\eeqq
The proof finishes by using Spitzer's Theorem (\ref{Spi}) and remarking that
\beq\label{alphacvg}
\frac{\log \alpha_{t}}{2\lambda t}\stackrel{t\rightarrow\infty}{\longrightarrow}1 .
\eeq
\end{proof}
{\noindent We finish this Subsection by stating and proving the following Large time asymptotic result
for the radial process of an Ornstein-Uhlenbeck process: }
\begin{theo}\label{radiaOUlarget}
The following convergence in law holds:
\beq\label{OUcvgradial}
\frac{\log R^{Z}_{t}}{t} \;  \overset{{(law)}}{\underset{t\rightarrow\infty}\longrightarrow} 0 .
\eeq
\end{theo}
\begin{proof}
From (\ref{RBMOU}), applying the scaling property of BM, we get
\beqq
\frac{\log|Z_{t}|}{t}=-\lambda+\frac{\log|B_{\alpha(t)}|}{t}
\stackrel{(law)}{=}-\lambda+\frac{\log(\sqrt{\alpha_{t}})+\log|B_{1}|}{t} .
\eeqq
Using (\ref{alphacvg}) and because $\lambda$ is a constant, we obtain (\ref{OUcvgradial}).
\end{proof}

\subsection{A complementary Large time asymptotics result}\label{subsectlarge}
Concerning the asymptotic behavior of the exit time from a cone with single boundary when $t\rightarrow\infty$,
we have the following:
\begin{prop}\label{proptlarge}
The asymptotic equivalence
\beqq
  2\lambda t \; P(T^{\theta(\lambda)}_{c}>t) \overset{{t\rightarrow \infty}}{\longrightarrow} \frac{4c}{\pi} ,
\eeqq
holds. It follows that, with $a,b>0$,
\beqq
2\lambda t \ P\left(a<\theta^{Z}_{t}<b\right)\overset{{t\rightarrow \infty}}{\longrightarrow} \frac{2}{\pi} (b-a).
\eeqq
\end{prop}
\begin{proof}
The first assertion follows from equation (\ref{Tlambda}) together with the analogous result for planar BM, that is
\beqq
  \left(\log t\right) \; P(T^{\theta}_{c}>t) \overset{{t\rightarrow \infty}}{\longrightarrow} \frac{4c}{\pi} ,
\eeqq
For the proof of the latter, see e.g. Proposition 2.5 in \cite{Vak11}.
Note that for this proof we could also invoke standard arguments, found e.g. in Pap-Yor \cite{PaY00} or a more recent proof
based on mod-convergence \cite{DKN11}.
The second convergence follows easily by remarking that
\beqq
2\lambda t \ P\left(a<|\theta^{Z}_{t}|<b\right)&=& 2\lambda t \ \left(P(T^{\theta}_{b}>\alpha_{t}) - P(T^{\theta}_{a}>\alpha_{t})\right) \nonumber \\
&\overset{{t\rightarrow \infty}}{\longrightarrow}& \frac{4}{\pi} (b-a),
\eeqq
and
\beqq
P\left(a<\theta^{Z}_{t}<b\right)=\frac{1}{2}P\left(a<|\theta^{Z}_{t}|<b\right).
\eeqq
\end{proof}

\subsection{Windings of complex-valued OU processes in $\left(\right.t,1\left.\right]$ for $t\rightarrow0$}\label{remt1}
We finish this Section by a study of complex-valued OU processes in a time interval.
Consider a 2-dimensional OU process $\left(\hat{Z}_{t},t\geq0\right)$ starting from 0 and we want to study its windings
in $\left(\right.t,1\left.\right]$ for $t\rightarrow0$. First, we remark that it doesn't visit again the origin but
it winds a.s. infinitely often around it. We denote $\left(\theta^{Z}_{(t,1)},0\leq t\leq1\right)$ its (well defined) continuous winding process
in the interval $\left(\right.t,1\left.\right]$, $t\leq1$.
We also denote by $\left(\hat{B}_{t},t\geq0\right)$ the planar BM starting from 0, which is associated to $\hat{Z}$.
\begin{prop} The following convergence in law holds:
\beqq
t \ \theta^{Z}_{(t,1)}\overset{{(law)}}{\underset{t\rightarrow0}\longrightarrow} C_{\lambda} .
\eeqq
\end{prop}
\begin{proof}
Changing variables $u=\alpha_{t}v$ and applying the scaling property of BM: $B_{\alpha_{t}v}\stackrel{(law)}{=}\sqrt{\alpha_{t}}\hat{B}_{v}$,
with obvious notation, identity (\ref{thetaBMOU}) yields
\beqq
\theta^{Z}_{(t,1)}=\theta^{B}_{(\alpha_{t},1)}=\mathrm{Im} \left(\int^{1}_{\alpha_{t}} \frac{dB_{u}}{B_{u}}\right)&\stackrel{(law)}{=}& \mathrm{Im} \left(\int^{1/\alpha_{t}}_{1} \frac{d\hat{B}_{v}}{\hat{B}_{v}}\right)= \theta^{\hat{B}}_{(1,1/\alpha_{t})}=\theta^{\hat{Z}}_{(1,1/t)} .
\eeqq
Hence, from Theorem \ref{SpitzerOU} we obtain
$$t \ \theta^{\hat{Z}}_{(1,1/t)}\overset{{(law)}}{\underset{t\rightarrow0}\longrightarrow} C_{\lambda},$$
which finishes the proof.
\end{proof}
\begin{rem}
For similar results concerning the windings of planar BM and (respectively of planar stable processes)
in $\left(\right.t,1\left.\right]$ for $t\rightarrow0$, see \cite{LeG92,ReY99} (respectively \cite{DoV12}). \\
Note that for the BM case, we can also invoke a time inversion argument (i.e.: $B_{u}=uB'_{1/u}$ where
$B'$ is another planar BM associated to an OU process $Z'$).
Hence, this argument could also be applied for the OU case studied here, i.e.
\beqq
\theta^{Z}_{(t,1)}=\theta^{B}_{(\alpha_{t},1)}&=&\mathrm{Im} \left(\int^{1}_{\alpha_{t}} \frac{dB_{u}}{B_{u}}\right)
=\mathrm{Im} \left(\int^{1}_{\alpha_{t}} \frac{d(uB'_{1/u})}{uB'_{1/u}}\right)
=\mathrm{Im} \left(\int^{1}_{\alpha_{t}} \frac{d(B'_{1/u})}{B'_{1/u}}\right) \nonumber \\
&=&\theta^{B'}_{(1,1/\alpha_{t})}=\theta^{Z'}_{(1,1/t)} ,
\eeqq
and we apply Theorem \ref{SpitzerOU} as before.
\end{rem}

\section{Limit Theorems for the exit time from a cone}\label{seclambda}

\subsection{Small and Big parameter asymptotics}\label{subseclambda}
We shall make use of the previously introduced notation for the first hitting times of a level $k>0$ for a real BM $\gamma$,
that is: $T^{\gamma}_{k}= \inf \left\{t\geq 0 : \gamma_{t}=k \right\}$ and
$T^{|\gamma|}_{k}= \inf \left\{t\geq 0 : |\gamma_{t}|=k \right\}$.
The following Proposition comes from \cite{Vak11} and we refer the reader therein for the proof and for further results.
\begin{prop}\label{OUlambda} For $z_{0}=1+i0$, the following convergence holds:
\beq\label{ETchatasymp}
    2\lambda \: E\left[T^{|\theta(\lambda)|}_{c}\right] - \log\left( 2\lambda \right) \stackrel{\lambda\rightarrow \infty}{\longrightarrow} E\left[\log\left(T^{|\theta|}_{c}\right)\right],
\eeq
with
\beqq
    E\left[\log\left(T^{|\theta|}_{c}\right)\right] = 2 \int^{\infty}_{0} \frac{dz}{\cosh \left(\frac{\pi z}{2}\right)} \log\left(\sinh\left(cz\right)\right) +\log\left(2\right) + c_{E},
\eeqq
where $c_{E}$ is Euler's constant. \\
For $c<\frac{\pi}{8}$, we also have the following convergence:
\beq\label{ETchatasymp02}
    \frac{1}{\lambda}\left(E\left[T^{|\theta(\lambda)|}_{c}\right] - E\left[\left(\sinh\left(\beta_{T^{|\gamma|}_{c}}\right)\right)^{2}\right]\right) \stackrel{\lambda\rightarrow 0}{\longrightarrow} - \frac{1}{3}E\left[\left(\sinh\left(\beta_{T^{|\gamma|}_{c}}\right)\right)^{4}\right].
\eeq
Equivalently,
\beq\label{ETchatderiv}
    \frac{d}{d\lambda} \Big|_{\lambda=0} E\left[T^{|\theta(\lambda)|}_{c}\right] = \underset{\lambda\rightarrow 0}{\lim} \left[\frac{1}{\lambda}\left(E\left[T^{|\theta(\lambda)|}_{c}\right] - E\left[T^{|\theta(0)|}_{c}\right]\right)\right] = - \frac{1}{3}E\left[\left(\sinh\left(\beta_{T^{|\gamma|}_{c}}\right)\right)^{4}\right].
\eeq
\end{prop}
\begin{rem}
We cannot get an analogue of (\ref{ETchatasymp}) for $E\left[T^{\theta(\lambda)}_{c}\right]$,
because the latter explodes, for every $c>0$.
Observe that the obvious analogs of formulae (\ref{ETchatasymp02}) and (\ref{ETchatderiv}) are not valid for $T^{\theta(\lambda)}_{c}$
for similar reasons.
\end{rem}

\subsection{Small and Big angle asymptotics}\label{subsecc}
In this Subsection, we study $T^{|\theta(\lambda)|}_{c}$ and $T^{\theta(\lambda)}_{c}$
for $c\rightarrow0$ and for $c\rightarrow\infty$ in the spirit of \cite{VaY12} (see also \cite{LGY86}).
Our main result is the following:
\begin{prop}\label{prop1}
$\left.a\right)$ For $c\rightarrow 0$, we have
\beqq
\frac{1}{c^{2}} \; T^{|\theta(\lambda)|}_{c} \overset{{(law)}}{\underset{c\rightarrow 0}\longrightarrow} T^{|\gamma|}_{1}.
\eeqq
$\left.b\right)$ For $c\rightarrow \infty$, we have
\beqq
\lambda \ \frac{T^{|\theta(\lambda)|}_{c}}{c} \overset{{(law)}}{\underset{c\rightarrow \infty}\longrightarrow} |\beta|_{T^{|\gamma|}_{1}}.
\eeqq
\end{prop}
\begin{proof} Both proofs are based on (\ref{Tabstheta}). \\
{\noindent$\left.a\right)$} It follows using the next elementary computation:
$$\frac{\log(1+2\lambda x)}{2\lambda}-x=\frac{1}{2\lambda}\int^{1+2\lambda x}_{0}\frac{dy}{y}
=\frac{1}{2\lambda}\int^{2\lambda x}_{0}\left(\frac{1}{1+a}-1\right)da.$$
Hence, taking $x=T^{|\theta|}_{c}$, recalling \eqref{Tabstheta} and invoking the fact that (Vakeroudis-Yor \cite{VaY12})
\beqq
\frac{1}{c^{2}} \; T^{|\theta|}_{c} \overset{{(law)}}{\underset{c\rightarrow 0}\longrightarrow} T^{|\gamma|}_{1},
\eeqq
we get
\beqq
\frac{1}{c^{2}}(T^{|\theta(\lambda)|}_{c}-T^{|\theta|}_{c})
=\frac{1}{2\lambda c^{2}} \int^{2\lambda T^{|\theta|}_{c}}_{0}\frac{-a}{1+a} \; da
=\int^{T^{|\theta|}_{c}/c^{2}}_{0}\frac{-2\lambda c^{2}}{1+2\lambda c^{2}} \; db .
\eeqq
Making $c\rightarrow 0$ in both sides, we get the announced result. \\ \\
{\noindent$\left.b\right)$} Obviously,
\beqq
2\lambda \ \frac{T^{|\theta(\lambda)|}_{c}}{c} &=& \frac{1}{c}\log \left(1+2\lambda T^{|\theta|}_{c}\right) \\
&=& \frac{1}{c}\log T^{|\theta|}_{c}+\frac{1}{c}\log \left(\frac{1}{T^{|\theta|}_{c}}+2\lambda \right).
\eeqq
The proof finishes by making $c\rightarrow \infty$ and using the result in \cite{VaY12}:
\beqq
\frac{1}{c}\log T^{|\theta|}_{c} \overset{{(law)}}{\underset{c\rightarrow \infty}\longrightarrow} |\beta|_{T^{|\gamma|}_{1}} .
\eeqq
\end{proof}
\begin{rem}
Comparing Proposition \ref{prop1} with Proposition 3.1 in \cite{VaY12}, we remark that the behavior
of the exit times from a cone of planar BM and of complex-valued OU processes is the same when $c\rightarrow0$
whereas it is different for $c\rightarrow\infty$.
\end{rem}
\subsection*{Generalizations}
Proposition \ref{prop1} has several variants. For instance
we define
$$T^{\theta(\lambda)}_{-b,a}=\inf \left\{t\geq 0 : \theta^{Z}_{t}\notin(b,a) \right\}, \ \ 0<a,b\leq \infty,$$
and
$$T^{\gamma}_{-d,c}=\inf\{ t:\gamma_{t}\notin(-d,c) \}, \ \ 0<c,d\leq \infty.$$
Hence, for $c\rightarrow 0$ or $c\rightarrow \infty$,
and $a,b$ fixed, we have
\begin{itemize}
\item $\frac{1}{c^{2}} \; T^{\theta(\lambda)}_{-bc,ac} \overset{{(law)}}{\underset{c\rightarrow 0}\longrightarrow} T^{\gamma}_{-b,a}$.
\item $ \lambda \ \frac{T^{\theta(\lambda)}_{-bc,ac}}{c} \overset{{(law)}}{\underset{c\rightarrow \infty}\longrightarrow} |\beta|_{T^{\gamma}_{-b,a}}$.
\end{itemize}
and with $b=\infty$, we get
\begin{cor}\label{coro}
$\left.a\right)$ For $c\rightarrow 0$, we have
\beqq
\frac{1}{c^{2}} \; T^{\theta(\lambda)}_{ac} \overset{{(law)}}{\underset{c\rightarrow 0}\longrightarrow} T^{\gamma}_{a}.
\eeqq
$\left.b\right)$For $c\rightarrow \infty$, we have
\beq\label{asympclarge}
\lambda \ \frac{T^{\theta(\lambda)}_{ac}}{c} \overset{{(law)}}{\underset{c\rightarrow \infty}\longrightarrow} |\beta|_{T^{\gamma}_{a}} \stackrel{(law)}{=}|C_{a}|,
\eeq
where $C_{a}$ is a Cauchy random variable.
\end{cor}
\begin{rem}\label{remSpi}  \emph{\textbf{(Yet another proof of Spitzer's Theorem for OU processes)}} \\
We remark that (\ref{asympclarge}) with $a=1$ yields another proof for the analogue of Spitzer's asymptotic Theorem
for OU processes (Theorem \ref{SpitzerOU}). Indeed, (\ref{asympclarge}) can be equivalently stated as:
\beq\label{asympP}
P\left( T^{\theta(\lambda)}_{c}<\frac{cx}{\lambda}\right)\overset{{(law)}}{\underset{c\rightarrow \infty}\longrightarrow}P\left(|C_{1}|<x\right).
\eeq
Invoking now the symmetry principle of Andr\'{e} \cite{And87,Gal08}, the LHS of (\ref{asympP}) is equal to
\beqq
P\left(\sup_{u\leq cx/\lambda} \theta^{Z}_{u} >c\right)
&=&P\left(\sup_{u\leq cx/\lambda} \theta^{B}_{\alpha(u)} >c\right)
=P\left(\sup_{u\leq cx/\lambda} \gamma_{H_{\alpha(u)}} >c\right) \nonumber \\
&=&P\left(|\gamma_{H_{\alpha(cx/\lambda)}}| >c\right)
=P\left(|\theta^{B}_{\alpha(cx/\lambda)}| >c\right)
=P\left(|\theta^{Z}_{cx/\lambda}| >c\right) \nonumber \\
&\stackrel{t=cx/\lambda}{=}&P\left(|\theta^{Z}_{t}| >\frac{\lambda t}{x}\right),
\eeqq
and (\ref{SpitzerOUcvg}) follows from (\ref{asympP}) for every $x>0$, by simply remarking that $|C_{1}|\overset{{(law)}}{=}|C_{1}|^{-1}$,
together with the fact that the symmetry principle yields again the following: for $k>0$,
\beqq
P\left(\theta^{Z}_{t}<k\right)&=&\frac{1}{2}P\left(|\theta^{Z}_{t}|<k\right), \\
P\left(C_{\lambda}<k\right)&=&\frac{1}{2}P\left(|C_{\lambda}|<k\right).
\eeqq
\end{rem}
\begin{rem}\label{remtheta}
Remark that the winding process of planar BM and that of complex-valued OU processes have the same behavior
when $c\rightarrow0$ limit, which is not the case when $c\rightarrow\infty$ (compare e.g. with \cite{VaY12}).
For some further results for the reciprocal of the exit time from a cone of planar Brownian motion $T^{|\theta|}_{c}$,
that is some infinite divisibility properties, see \cite{VaY12b}.
\end{rem}
\begin{rem}
The interested reader can also compare the results for the exit times from a cone
with the analogues of processes with jumps (stable processes) in \cite{DoV12}.
\end{rem}

\section{Small and Big windings of Ornstein-Uhlenbeck processes}\label{secbigsmall}

\subsection{Small and Big windings}\label{subsecbigsmall}
As for planar BM (see e.g. \cite{PiY84,PiY86,LGY86}), it is natural to continue
the study of the windings of complex-valued OU processes by decomposing the winding process in
"small" and "big" windings. To that direction, because of the positive recurrence of OU processes,
we expect a significantly different asymptotic behavior (when $t\rightarrow\infty$) of these two components
comparing to that of BM, which is null recurrent.

Following e.g. \cite{PiY86}, we consider $\mathbb{C}$ the whole complex domain where $Z$ a.s. "lives" and we decompose it in
$D_{+}$ (the big domain) and $D_{-}$ (the small domain) the open sets outside and inside the unit circle
(hence: $D_{+}+D_{-}=\mathbb{C} \setminus \left\{z:|z|=1\right\}$),
with the sign + and - standing for big and small respectively (inspired by the sign of $\log|z|$, with $z$ in the whole domain). We define
\beq\label{Zpm}
 \theta^{Z}_{\pm}(t) = \int^{t}_{0} \: 1(Z(s) \in D_{\pm}) \: d\theta^{Z}_{s},
\eeq
where $1(A)$ is the indicator of $A$.
The process $\theta^{Z}_{+}$ is the process of big windings and $\theta^{Z}_{-}$ is the process of small windings,
both associated to $Z$.
The Lebesgue measure of the time spent by $Z$ on the unit circle is a.s. 0, thus
\beq\label{Zdecomposition}
 \theta^{Z} = \theta^{Z}_{+}+\theta^{Z}_{-}.
\eeq
Recall that, as mentioned in Subsection \ref{largetimeasymp},
the (well-defined) winding process $\vartheta_{t}$
of a planar diffusion starting from a point different from the origin was studied by
Friedman and Pinsky in \cite{FrP73,FrP73a}, and they showed that, when $t\rightarrow\infty$,
$\vartheta_{t}/t$ exists a.s. under some assumptions implying that the process winds asymptotically around a point.

A first remark is that, similar to planar BM, the winding process $\theta$ is switching between long time periods, when $Z$ is far away from the origin in $D_{+}$ and $\theta$ changes very slowly (but significantly) because of $\theta_{+}$, and small time periods, when $Z$ is in $D_{-}$ approaching 0 and $\theta$ changes very rapidly because of $\theta_{-}$.
It follows that, contrary to planar BM where the very big windings and very small windings count for the asymptotic behavior (as $t\rightarrow\infty$) of the total winding,
for OU processes only the very small windings contribute. We also note that, the windings for a very large class of 2-dimensional random walks, behave rather more like $\theta_{+}$ than $\theta$ (see e.g. \cite{Bel86,Bel89,BeF91,BeW94,Shi98}).
\\ \\
First, we extend Theorem 1 (iii) in Bertoin and Werner \cite{BeW94}.
\begin{prop}\label{propwindings}
We consider $f$ a complex-valued bounded Borel function with compact support on the whole complex domain $\mathbb{C}$. Then,
with $z\in\mathbb{C}$ (equivalently $z=x+iy$), we have
\beq\label{intf}
\frac{1}{t}\int^{t}_{0}ds \ f(Z_{s}) \overset{{a.s.}}{\underset{t\rightarrow \infty}\longrightarrow}
\frac{\lambda}{\pi}\int_{\mathbb{R}^{2}} dx \; dy \;  e^{-\lambda(x^{2}+y^{2})} f(z) .
\eeq
\end{prop}
\begin{proof}
We start by noting that, for fixed $s$, $Z_{s}$ is bivariate normally distributed where each component has mean 0 and variance $\exp(-2\lambda s)\alpha_{s}$, where we also recall that:
$$\alpha_{s}=\frac{1}{2\lambda}\left(e^{2\lambda s}-1\right).$$
Hence, the variance converges to $1/(2\lambda)$ as $s\rightarrow\infty$, and we obtain the invariant probability measure
of $(Z_{t}, t\geq0)$, that is: $$\frac{\lambda}{\pi} \; e^{-\lambda |z|^{2}} dx \; dy .$$
Invoking the Ergodic Theorem, we obtain
\beqq
\frac{1}{t}\int^{t}_{0}ds \ f\left(Z_{s}\right) \overset{{a.s.}}{\underset{t\rightarrow \infty}\longrightarrow}
\int_{\mathbb{R}^{2}} dx \; dy \; \frac{\lambda}{\pi} \; e^{-\lambda |z|^{2}} f(z) ,
\eeqq
which is precisely (\ref{intf}).
\end{proof}
We consider now, without loss of generality, that $D_{+}$ and $D_{-}$ are such that
$|Z_{\cdot}|\in(1,+\infty)$ and $|Z_{\cdot}|\in(0,1)$ respectively.
Hence, using (\ref{thetaBMOU}), we may write
\beq\label{Zp}
 \theta^{Z}_{+}(t) &=& \int^{t}_{0} \: 1(|Z_{s}| \geq1) \: \mathrm{Im}\left(\frac{dZ_{s}}{Z_{s}}\right)
 =\int^{t}_{0} \: 1(|Z_{s}| \geq1) \: \mathrm{Im}\left(\frac{dB_{\alpha(s)}}{B_{\alpha(s)}}\right) \nonumber \\
 &=& \int^{\alpha(t)}_{0} \: 1(|Z_{\alpha^{-1}(u)}| \geq1) \: d\theta^{B}_{u} ,
\eeq
where, for the latter, we have changed the variables $u=\alpha(s)$.
Similarly,
\beqq
 \theta^{Z}_{-}(t) = \int^{t}_{0} \: 1(|Z_{\alpha^{-1}(u)}| \leq1) \: d\theta^{B}_{u} \ .
\eeqq
\begin{theo}\label{theosbwindings}
The following convergence in law holds:
\beq\label{triplet}
\frac{1}{t} \; \theta^{Z}_{+}(t)\overset{{(P)}}{\underset{t\rightarrow \infty}\longrightarrow}0 ,
\eeq
while
\beq\label{triplet2}
\frac{1}{t} \; \theta^{Z}_{-}(t)\overset{{(law)}}{\underset{t\rightarrow \infty}\longrightarrow}C_{\lambda} .
\eeq
\end{theo}
\begin{rem}
Theorem \ref{theosbwindings} essentially means that the big windings of complex-valued Ornstein-Uhlenbeck processes,
do not contribute to the total windings at the limit $t\rightarrow\infty$.
Hence, it is only the small windings that is taken into account at the large time limit, which seems natural
if we recall that OU processes are characterized by a force "pulling" them back to their origin,
thus they are positive recurrent.
\end{rem}
\begin{proof}
With $R^{Z}=(R^{Z}_{t}, t\geq0)=(|Z_{t}|, t\geq0)$, we define
(see also Section 2 in Bertoin and Werner \cite{BeW94} where a slightly different notation is used,
and \cite{MeY82,PiY86}): for every $\varepsilon>0$,
\beq\label{thetaepsilon}
\theta^{Z}_{\varepsilon}(e^{t})=\int^{t}_{0} 1_{\left(R(s)>\varepsilon\right)} d\theta^{Z}_{s} \ , \ \ t\geq1.
\eeq
Moreover, with $\varepsilon=0$, Spitzer's Theorem for OU processes (Theorem \ref{SpitzerOU}) yields
\beq\label{SpitzerOUcvgprime}
\frac{\theta^{Z}_{0}(e^{t})}{t}=\frac{\theta^{Z}(t)}{t} \;  \overset{{(law)}}{\underset{t\rightarrow\infty}\longrightarrow} C_{\lambda} .
\eeq
We will study now separately $\theta^{Z}_{+}$ and $\theta^{Z}_{-}$.
Note that we could use Proposition \ref{propwindings} in the spirit of Kallianpur-Robbins law
(we address the interested reader to e.g. Pitman-Yor \cite{PiY86}, or \cite{KaR53} for the original article).
However, we proceed to the following straightforward computations. \\ \\
{\noindent $\left.i\right)$} We start by equation (\ref{Zp}).
Using now (\ref{OUeB1}) and (\ref{alpha}), we have:
\beqq
 \theta^{Z}_{+}(t)&=&\int^{\alpha(t)}_{0} \: 1(e^{-\lambda \alpha^{-1}(u)}|B_{u}| \geq1) \: d\theta^{B}_{u}
 = \int^{\alpha(t)}_{0} \: 1(-\lambda \alpha^{-1}(u)+\log|B_{u}| \geq0) \: d\theta^{B}_{u} \\
 &=& \int^{\alpha(t)}_{0} \: 1\left(\log|B_{u}| \geq\frac{1}{2} \log(1+2\lambda u)\right) \: d\theta^{B}_{u} .
\eeqq
The skew-product representation (\ref{skew-product}) of the planar Brownian motion $B$ yields that (we also recall that
$A_{u}=A_{u}(\beta)=\int^{u}_{0} \exp(2\beta_{s})ds=H^{-1}_{u}$)
\beqq
\theta^{Z}_{+}(t)&=&\int^{\alpha(t)}_{0} \: 1\left(\beta_{H(u)} \geq\frac{1}{2} \log(1+2\lambda A_{H(u)})\right) \: d\gamma_{H(u)} \\
&\stackrel{v=H(u)}{=}&\int^{H_{\alpha(t)}}_{0} \: 1\left(\beta_{v} \geq\frac{1}{2}\log(1+2\lambda A_{v})\right) \: d\gamma_{v} .
\eeqq
On the one hand, with $\hat{\beta}$ and $\hat{\gamma}$ denoting two other real BMs starting from 0, independent from each other, such that: for every $t$, $\hat{\beta}_{w}=(\lambda t)^{-1}\beta_{\lambda^{2}t^{2}w}$ and $\hat{\gamma}_{w}=(\lambda t)^{-1}\gamma_{\lambda^{2}t^{2}w}$, and changing the variables $v=\lambda^{2}t^{2}w$, we obtain
\beq\label{eqproof}
&& \frac{1}{t} \int^{H_{\alpha(t)}}_{0} 1\left(\beta_{v} \geq\frac{1}{2}\log(1+2\lambda A_{v})\right) d\gamma_{v} \nonumber \\
&& \ \ \ \ \ \ \ \ \ \ \ \ \ \ \ \ \ \ \ \
=\lambda \int^{\frac{1}{\lambda^{2}t^{2}}H_{\alpha(t)}}_{0} 1\left(\hat{\beta}_{w} \geq\frac{1}{2\lambda t} \log(1+2\lambda A_{\lambda^{2}t^{2}w})\right) d\hat{\gamma}_{w} .
\eeq
Moreover,
\beqq
\frac{1}{t^{2}} H_{\alpha(t)}=\frac{1}{t^{2}} H_{\left(\frac{\exp(2\lambda t)-1}{2\lambda}\right)} ,
\eeqq
and recalling that (see e.g. \cite{LGY86,ReY99})
\beqq
\frac{4}{(\log u)^{2}} \; H_{u} \overset{{(law)}}{\underset{u\rightarrow \infty}\longrightarrow} T^{\beta}_{1}=\inf\{t:\beta_{t}=1\}=\frac{1}{N^{2}}, \ \ \mathrm{with} \ \ N\thicksim\mathcal{N}(0,1)
\eeqq
we get
\beq\label{Hbigcvg}
\frac{1}{\lambda^{2}t^{2}} H_{\left(\frac{\exp(2\lambda t)-1}{2\lambda}\right)}\overset{{(law)}}{\underset{t\rightarrow \infty}\longrightarrow} T^{\beta}_{1} .
\eeq
On the other hand, changing the variables $s=\lambda^{2}t^{2}u$,
\beq\label{alphainv}
&& \frac{1}{2\lambda t} \log(1+2\lambda A_{\lambda^{2}t^{2}w})=
\frac{1}{2\lambda t} \log\left(1+2\lambda \int^{\lambda^{2}t^{2}w}_{0}e^{2\beta_{s}}ds\right) \nonumber \\
&=&\frac{1}{2\lambda t} \log\left(1+2\lambda^{3} t^{2} \int^{w}_{0}e^{2\lambda t\hat{\beta}_{u}}du\right) \nonumber \\
&=& \frac{\log(2\lambda^{3} t^{2})}{2\lambda t}+ \frac{1}{2\lambda t} \log\left(1+ \left(2\lambda^{3} t^{2} \int^{w}_{0}e^{2\lambda t\hat{\beta}_{u}}du\right)^{-1}\right)
+ \log\left(\int^{w}_{0}e^{2\lambda t\hat{\beta}_{u}}du\right)^{1/(2\lambda t)} \nonumber \\
&\overset{{(P)}}{\underset{t\rightarrow \infty}\longrightarrow}&
\log\left(\sup_{u\leq w} e^{\hat{\beta}_{u}}\right)=\sup_{u\leq w} \hat{\beta}_{u} ,
\eeq
where the latter follows by invoking again the convergence of the $p$-norm to the $\infty$-norm,
as $p\rightarrow\infty$.
Convergence (\ref{alphainv}), together with (\ref{eqproof}) and (\ref{Hbigcvg}), yields that
\beqq
\frac{\theta^{Z}_{+}(t)}{t}\overset{{(law)}}{\underset{t\rightarrow \infty}\longrightarrow}
\int^{T^{\beta}_{1}}_{0} 1\left(\hat{\beta}_{w} \geq \sup_{u\leq w} \hat{\beta}_{u}\right) d\hat{\gamma}_{w}=0 ,
\eeqq
hence, it also converges to 0 in Probability. \\ \\
{\noindent $\left.ii\right)$}
Concerning the small windings process $\theta^{Z}_{-}$,
the decomposition in small and big windings (\ref{Zdecomposition})
together with Spitzer's Theorem for OU processes (Theorem \ref{SpitzerOU}-or equivalently (\ref{SpitzerOUcvgprime}) )
and convergence in Probability (\ref{triplet}) for the big windings, yield (\ref{triplet2}).
\\ \\
We note that for part $\left.ii\right)$ of the proof, we could also mimic the proof for the Brownian motion case
(see e.g. \cite{PiY86} and in particular Lemma 3.1 and Theorem 4.1 therein),
invoking Williams "pinching method". This method was introduced in \cite{Wil74} and further investigated in \cite{MeY82}
(for other variations, see also \cite{Dur82,Dur84}).
\end{proof}
\begin{rem}
From (\ref{thetaepsilon}), using the skew-product representation and the Ergodic Theorem (as in the proof of Theorem 1 (iii) in \cite{BeW94}),
and recalling that $(1/2)1_{(u\geq0)}e^{-\lambda u} du$ is the invariant probability measure of $R^{2}$, we get
\beqq
\frac{\theta^{Z}_{\varepsilon}(e^{t})}{\sqrt{t}} \overset{{(law)}}{\underset{t\rightarrow \infty}\longrightarrow} k_{\varepsilon} \mathcal{N} \ ,
\eeqq
where $ k^{2}_{\varepsilon}=\int^{\infty}_{\varepsilon^{2}} u^{-1}e^{-\lambda u}du$ and $\mathcal{N}\thicksim N(0,1)$.
\end{rem}
\begin{rem}
We finish this Subsection by remarking that, as already mentioned in Bertoin-Werner \cite{BeW94} (see the Introduction therein),
contrary to the planar Brownian motion, this method does not seem to apply to the windings
of a complex-valued Ornstein-Uhlenbeck process about several points.
\end{rem}

\subsection{Very Big Windings}\label{verybig}
Theorem \ref{theosbwindings} (and in particular part $\left.i\right)$) is corresponding to the discussion already made in Bertoin-Werner \cite{BeW94}
where they introduced the $\nu$-big (respectively $\nu$-small) windings of planar BM (we use a slightly modified notation
convenient for the needs of the present work), i.e.
\beqq
\theta^{B,\nu}_{t}&=& \int^{t}_{1} \: 1(|B_{s}| \geq s^{\nu}) \: d\theta^{B}_{s} \ , \ \ t\geq 1 ; \\
\theta^{B,-\nu}_{t}&=& \int^{t}_{1} \: 1(|B_{s}| \leq s^{-\nu}) \: d\theta^{B}_{s} \ , \ \ t\geq 1 ,
\eeqq
and saying that the case $\nu=1/2$ is a critical case which corresponds to the so-called
very big windings $\theta^{B,1/2}$ (see also Le Gall-Yor \cite{LGY90}). \\
Indeed, repeating the arguments of part $\left.i\right)$ in the proof of Theorem \ref{theosbwindings} with some modifications
(e.g. in the equation corresponding to (\ref{eqproof}), change the variables $u=(\log t)^{2}w$), we get
\beqq
\theta^{B,\nu}_{t}\overset{{(law)}}{\underset{t\rightarrow \infty}\longrightarrow}\int^{T^{\beta}_{1}}_{0} 1\left(\beta_{v}\geq0\right) d\gamma_{v}
\ \Longleftrightarrow \ \nu<1/2.
\eeqq
We turn now our study to the $\nu$-big (respectively $\nu$-small) windings of complex-valued OU processes
\beqq
\theta^{Z,\nu}_{t}&=& \int^{\alpha(t)}_{1} \: 1(|Z_{s}| \geq s^{\nu}) \: d\theta^{B}_{s} \ , \ \ t\geq 1 ; \\
\theta^{Z,-\nu}_{t}&=& \int^{\alpha(t)}_{1} \: 1(|Z_{s}| \leq s^{-\nu}) \: d\theta^{B}_{s} \ , \ \ t\geq 1 .
\eeqq
\begin{prop} The following convergence in law holds:
\beq\label{thetaverybig}
\theta^{Z,\nu}_{t}&\overset{{(law)}}{\underset{t\rightarrow \infty}\longrightarrow}&
\int^{T^{\hat{\beta}}_{1}}_{0} 1\left(\hat{\beta}_{v}\geq (1+2\nu) \sup_{u\leq v} \hat{\beta}_{u}\right) d\gamma_{v} ,
\eeq
which is not degenerate if and only if
$1+2\nu<1 \ \Longleftrightarrow \ \nu<0$,
and
\beq\label{thetaverysmall}
\theta^{Z,-\nu}_{t}&\overset{{(law)}}{\underset{t\rightarrow \infty}\longrightarrow}&
\int^{T^{\hat{\beta}}_{1}}_{0} 1\left(\hat{\beta}_{v}\leq (1-2\nu) \sup_{u\leq v} \hat{\beta}_{u}\right) d\gamma_{v}
\eeq
which is not degenerate if and only if
$1-2\nu<1 \ \Longleftrightarrow \ \nu>0$.
\end{prop}
\begin{proof}
The slightly modified arguments above in the proof of Theorem \ref{theosbwindings} yield that
\beqq
\theta^{Z,\nu}_{t}=\int^{H_{\alpha(t)}}_{0} \: 1\left(\beta_{v} \geq\frac{1}{2}\log(1+2\lambda A_{v})+\nu \log A_{v}\right) \: d\gamma_{v} ,
\eeqq
and
\beqq
\frac{1}{2t}\log(1+2\lambda A_{v})+\frac{\nu}{t}\log A_{v} &\stackrel{v=t^{2}w}{=}&
\frac{1}{2t}\log(1+2\lambda A_{t^{2}w})+\frac{\nu}{t}\log A_{t^{2}w} \nonumber \\
&\overset{{(P)}}{\underset{t\rightarrow \infty}\longrightarrow}& (1+2\nu) \sup_{u\leq w} \hat{\beta}_{u} .
\eeqq
hence we get \eqref{thetaverybig}. Similarly, we obtain \eqref{thetaverysmall}.
\end{proof}

\section{Windings of Ornstein-Uhlenbeck processes driven by a Stable process (OUSP)}\label{OUL}

\subsection{Preliminaries on L\'evy and Stable processes}\label{LevyIntro}
For some basic properties of L\'{e}vy processes and Stable processes
we refer to e.g. \cite{Ber96} or \cite{Kyp06}.

Coming from Lamperti \cite{Lam72}, a Markov process $J$ taking values in $\mathbb{R}^{d}$, $d\geq2$
is called \textit{isotropic} or $O(d)$-\textit{invariant} ($O(d)$ is the group of orthogonal transformations
on $\mathbb{R}^{d}$) if its transition satisfies
\beqq
P_{t}(\phi(x),\phi(\mathcal{B}))=P_{t}(x,\mathcal{B}),
\eeqq
for any $\phi\in O(d)$, $x\in \mathbb{R}^{d}$ and Borel subset $\mathcal{B}\subset\mathbb{R}^{d}$. \\
Moreover, $J$ is said to be $\alpha$-\textit{self-similar} if, for $\alpha>0$,
\beqq
P_{\psi t}(x,\mathcal{B})=P_{t}(\psi^{-\alpha}x,\psi^{-\alpha}\mathcal{B}),
\eeqq
for any $\psi>0$, $x\in \mathbb{R}^{d}$ and $\mathcal{B}\subset\mathbb{R}^{d}$.

We turn now our interest to the 2-dimensional case $(d=2)$.
We denote by $(\tilde{U}_{t},t\geq0)$ a standard isotropic stable process of index $\alpha\in(0,2)$
taking values in the complex plane and starting from $u_{0}+i0,u_{0}>0$.
Without loss of generality (it follows easily by a scaling argument), from now on we may assume that $u_{0}=1$.
Some basic properties of $\tilde{U}$ are the following (see e.g. \cite{Ber96,Kyp06}): it has stationary independent increments,
its sample paths are right continuous and has left limits (cadlag) and, with $\langle\cdot,\cdot\rangle$ standing for the Euclidean inner product, $E\left[\exp\left(i\langle\lambda,\tilde{U}_{t}\rangle\right)\right]=\exp\left(-t |\lambda|^{\alpha}\right)$, for all $t\geq0$ and $\lambda\in\mathbb{C}$.
$\tilde{U}$ is transient, $\lim_{t\rightarrow\infty} |\tilde{U}_{t}|=\infty$ a.s.
and it a.s. never visits single points.
Note that for $\alpha=2$, we are in the Brownian motion case.

We also introduce the following processes: $Q=(Q_{t},t\geq0)$ denotes a planar Brownian motion starting from $1+i0$ and
$S=(S(t),t\geq0)$ stands for an independent stable subordinator with index $\alpha/2$ starting from 0, where $\alpha\in(0,2)$, i.e.
\beqq
E\left[\exp\left(-\mu S(t)\right)\right]=\exp\left(-t \mu^{\alpha/2}\right),
\eeqq
for all $t\geq0$ and $\mu\geq0$. It follows that the subordinated planar Brownian motion $\tilde{U}_{\cdot}=Q_{2S(\cdot)}$
is a standard isotropic stable process of index $\alpha$.
The L\'{e}vy measure of $S$ is
$$ \frac{\alpha}{2\Gamma(1-\alpha/2)} \ s^{-1-\alpha/2} 1_{\{s>0\}} ds \ .$$
and it follows that, the L\'{e}vy measure $\nu$ of $\tilde{U}$ is (see e.g. \cite{BeW96})
\beqq
\nu(dx)&=& \frac{\alpha}{2\Gamma(1-\alpha/2)} \int^{\infty}_{0} s^{-1-\alpha/2} P\left(Q_{2s}-1 \ \in \ dx\right) ds \nonumber \\
&=& \frac{\alpha}{8\pi\Gamma(1-\alpha/2)} \left(\int^{\infty}_{0} s^{-2-\alpha/2} \exp\left(-|x|^{2}/(4s)\right) \ ds\right) dx \nonumber \\
&=&\frac{\alpha \ 2^{-1+\alpha/2} \Gamma(1+\alpha/2)}{\pi\Gamma(1-\alpha/2)} \ |x|^{-2-\alpha} dx .
\eeqq
The windings of Stable processes have already been studied and we refer the interested reader to
Bertoin-Werner \cite{BeW96}, Doney-Vakeroudis \cite{DoV12} and the references therein.

\subsection{Windings of planar OU processes driven by a BDL process}\label{subsecLevy}
We turn now our study to the windings of complex-valued Ornstein-Uhlenbeck processes driven by a Stable process (OUSP).
We consider
\beq\label{OULeq}
    V_{t} = v_{0} + U_{\lambda t} - \lambda \int^{t}_{0} V_{s} ds,
\eeq
with $\left(U_{t},t\geq0\right)$ denoting the Background 2-dimensional time homogeneous driving L\'evy (Stable in our case) process (BDLP), starting from 0, a terminology initially introduced in \cite{BNS01}, $v_{0}\in \mathbb{C}^{\ast}$ and $\lambda \geq 0$ (for more details about BDLP, see also \cite{Spi09,Ona09} and the references therein). Note that, following \cite{BNS01} p. 175, the SDE satisfied by $V$
is written in the form (\ref{OULeq}), which follows after a simple change of variables, in order to obtain a stationary solution. \\
We also have the following representation:
\beqq
    V_{t} &=& e^{-\lambda t} \left( v_{0} + \int^{\lambda t}_{0} e^{s} dU_{s} \right),
\eeqq
which is equivalent to (\ref{OULeq}) by using e.g. It\^{o}'s formula. \\
Without loss of generality, we may suppose: $v_{0}=1+i0$.
Moreover, writing now $U$ as a subordinated planar BM, i.e.: $Q_{2S(t)}$, we obtain
\beqq
V_{t}&=&e^{-\lambda t} \left( 1 + \int^{\lambda t}_{0} e^{s} dQ_{2S(s)} \right) .
\eeqq
We use now: $\left(V_{t}=V^{(1)}_{t}+iV^{(2)}_{t}; t\geq 0\right)$ and
$\left(U_{t}=U^{(1)}_{t}+iU^{(2)}_{t}; t\geq 0\right)$, where
$V^{(1)},V^{(2)}$ are two independent 1-dimensional OU processes starting respectively from 1 and 0,
and $U^{(1)},U^{(2)}$
are two independent 1-dimensional Stable processes (with the same index of stability $\alpha$) starting both from 0.
As $V$ starts from a point different from 0, following \cite{BeW96} or \cite{DoV12},
we can consider a path on a finite time interval $[0,t]$ and "fill in" the gaps with line segments.
In that way, we obtain the curve of a continuous function
$f:[0,1]\rightarrow\mathbb{C}$ with $f(0)=1$ and since 0 is polar and $V$ has no jumps across 0 a.s.,
its winding process $\theta^{V}=\left(\theta^{V}_{t}, t\geq0\right)$ is well defined.
\begin{prop}\label{OULwindingsSDE}
The winding and the radial process of a complex-valued OU process $V$ driven by a Stable process satisfy respectively the following SDEs:
\beq
\theta^{V}_{t} &=& \lambda^{1/\alpha}\int^{t}_{0}\frac{V^{(1)}_{s}dU^{(2)}_{s}-V^{(2)}_{s}dU^{(1)}_{s}}{|V_{s}|^{2}} , \label{thetaBMOULSDE2} \\
\log R^{V}_{t} &=& -\lambda t+\lambda^{1/\alpha}\int^{t}_{0}\frac{V^{(1)}_{s}dU^{(1)}_{s}+V^{(2)}_{s}dU^{(2)}_{s}}{|V_{s}|^{2}} . \label{RBMOULSDE2}
\eeq
\end{prop}
\begin{proof}
We start by writing (\ref{OULeq}) in differential form, i.e.
\beqq
dV_{t} =  dU_{\lambda t} - \lambda V_{t} dt, \ \ \ \ V_{0}=v_{0}=1+i0 .
\eeqq
Hence,
\beqq
\mathrm{Im}\left(\frac{dV_{t}}{V_{t}}\right)&=&\mathrm{Im}\left(\frac{dU_{\lambda t} - \lambda V_{t} dt}{V_{t}}\right)
=\mathrm{Im}\left(\frac{dU_{\lambda t}}{V_{t}}\right)
=\mathrm{Im}\left(\frac{d\left(U^{(1)}_{\lambda t}+iU^{(2)}_{\lambda t}\right)}{V^{(1)}_{t}+iV^{(2)}_{t}}\right) \\
&=&\frac{-V^{(2)}_{t}dU^{(1)}_{\lambda t}+V^{(1)}_{t}dU^{(2)}_{\lambda t}}{|V_{t}|^{2}} ,
\eeqq
which writes
\beqq
\theta^{V}_{t} &=& \lambda^{1/\alpha}\int^{t}_{0}\frac{V^{(1)}_{s}dU^{(2)}_{s}-V^{(2)}_{s}dU^{(1)}_{s}}{|V_{s}|^{2}} ,
\eeqq
and equation (\ref{thetaBMOULSDE2}) follows by applying the stability property:
$U^{(j)}_{\lambda t}\stackrel{(law)}{=}\lambda^{1/\alpha} U^{(j)}_{t}$, $j=1,2$.
Similar computations for the radial process $\left(R^{V}_{t}=|V_{t}|, t\geq0\right)$, yield
\beqq
\log R^{V}_{t} &=& -\lambda t+\int^{t}_{0}\frac{V^{(1)}_{s}}{|V_{s}|^{2}} \ dU^{(1)}_{\lambda s}+
\int^{t}_{0}\frac{V^{(2)}_{s}}{|V_{s}|^{2}} \ dU^{(2)}_{\lambda s}
\eeqq
thus \eqref{RBMOULSDE2}.
\end{proof}

\subsection{Windings of planar OU processes driven by a Stable process}\label{secdevelop}
In this last Subsection, we will investigate the case of the complex-valued OU process
\beqq
    V_{t} = v_{0} + J_{t} - \lambda \int^{t}_{0} V_{s} ds,
\eeqq
where $(J_{t})_{t\geq0}$ is an $\alpha$-stable process with $\alpha\in(0,2]$.
We also introduce the clock:
\beqq
H^{J}_{t}\equiv\int^{t}_{0}\frac{ds}{\left|J_{s}\right|^{\alpha}} \ ,
\eeqq
having as an inverse:
\beq\label{HAclock}
(H^{J})^{-1}_{u}\equiv A^{J}_{u}\equiv \inf\{t\geq0: H^{J}_{t}>u\}=\int^{u}_{0} \exp\{\alpha \xi_{s}\} \ ds \ .
\eeq
Following \cite{BeW96}, we may get the Lamperti correspondence for stable processes (the analogue of the skew product representation for planar BM).
Indeed, following \cite{GVA86} and using Lamperti's relation (see e.g. \cite{ReY99}),
there exist two real-valued L\'{e}vy processes $(\xi_{u},u\geq0)$ and $(\rho_{u},u\geq0)$, where the first one is non-symmetric
whereas the second one is symmetric, both starting from 0, such that:
\beq\label{skew-productstable}
\log\left|J_{t}\right|+i\theta^{J}_{t}
=\left(\xi_{u}+i\rho_{u}\right)
\Bigm|_{u=H^{J}_{t}=\int^{t}_{0}\frac{ds}{\left|J_{s}\right|^{\alpha}}} \ .
\eeq
Note here that, contrary to the BM case,
$|J|$ and $J_{A^{J}_{\cdot}}/|J_{A^{J}_{\cdot}}|$ are not independent as,
roughly speaking, they jump at the same times
(see \cite{BeW96,DoV12} and the references therein for further discussion).
Using \eqref{HAclock}, from (\ref{skew-productstable}) we get
\beq\label{skew-productstable2}
\left\{
  \begin{array}{ll}
    \left|J_{t}\right|=\exp\left(\xi(H^{J}_{t})\right)\Leftrightarrow\left|J_{A^{J}_{t}}\right|
    =\exp\left(\xi_{t}\right), & \hbox{(extension of Lamperti's identity)} \\
    \theta^{J}_{t}=\rho(H^{J}_{t})\Leftrightarrow \theta\left(A^{J}_{t}\right)=\rho(t) \ .
  \end{array}
\right.
\eeq
We also define the random times $T^{|\theta^{J}|}_{c}\equiv\inf\{ t:|\theta^{J}_{t}|\geq c \}$ and
$T^{|\rho|}_{c}\equiv\inf\{ t:|\rho_{t}|\geq c \}$, with $c>0$,
and the "generalized" skew-product representation (\ref{skew-productstable}) (or (\ref{skew-productstable2})) writes:
\beqq
T^{|\theta^{J}|}_{c}=(H^{J})^{-1}_{u}\Bigm|_{u=T^{|\rho|}_{c}}=\int^{T^{|\rho|}_{c}}_{0}ds\exp(\alpha\xi_{s})\equiv A^{J}_{T^{|\rho|}_{c}} \ .
\eeqq
\begin{prop} The following relation holds:
\beq\label{relstable}
\theta^{V}_{t}=\theta^{\widetilde{J}}_{\widehat{\alpha}(t)},
\eeq
where
$$\widehat{\alpha}(t)=\int^{t}_{0}e^{\alpha \lambda s}ds=\frac{e^{\alpha \lambda t}-1}{\alpha \lambda}
\ \Longleftrightarrow \ \widehat{\alpha}^{-1}(t)=\frac{1}{\alpha \lambda} \ \log(1+\alpha \lambda t).$$
Hence, the hitting time
$\tau^{|\theta^{V}|}_{c}=\inf\{t\geq0, |\theta^{V}_{t}|=c\}$
satisfies
\beq\label{reltransforms}
\tau^{|\theta^{V}|}_{c}=\frac{1}{\alpha \lambda} \ \log(1+\alpha \lambda \tau^{|\theta^{\widetilde{J}}|}_{c}),
\eeq
or equivalently
\beq\label{transforms}
\mathbb{E}\left[e^{-u\tau^{|\theta^{V}|}_{c}}\right]=\mathbb{E}\left[(1+\alpha \lambda \tau^{|\theta^{\widetilde{J}}|}_{c})^{-u/(\alpha \lambda)}\right] .
\eeq
\end{prop}
\begin{proof}
First, we use Dubins-Schwartz Theorem which extends to the case of $\alpha$-stable processes
(see e.g. \cite{Kal06,KaS02}), meaning that there exists an independent $\alpha$-stable process $\widetilde{J}$ starting from $v_{0}$ such that
$$v_{0}+\int^{t}_{0}e^{\lambda s}dJ_{s}=\widetilde{J}_{\widehat{\alpha}(t)}.$$
Similar computations as in the complex-valued OU driven by a BM case (see e.g. \cite{Vak11}) yield \eqref{relstable}.
Using the latter we get (with obvious notation)
\beqq
\tau^{|\theta^{V}|}_{c}=\inf\{t\geq0, |\theta^{\widetilde{J}}_{\widehat{\alpha}(t)}|=c\}=\widehat{\alpha}^{-1}(\tau^{|\theta^{\widetilde{J}}|}_{c}),
\eeqq
thus \eqref{reltransforms}. Finally, \eqref{transforms} follows from \eqref{reltransforms} by taking the Laplace transform in both sides.
\end{proof}
{\noindent Mimicking the study of complex-valued OU processes driven by BM, we can obtain similar asymptotic results by invoking the asymptotics of stable processes from \cite{BeW96,DoV12}. In particular, using the "generalized" skew product representation \eqref{skew-productstable2} together with Theorems 4.4 and 3.2 from \cite{DoV12} respectively, we get the following small and big time asymptotics. Note that both results below refer to convergence in distribution on $D(\left[\right.0,\infty\left.\right),\mathbb{R})$ endowed with the Skorohod topology.}
\begin{theo}
\begin{enumerate}[(i)]
  \item The family of processes $(t^{-1/\alpha}\theta^{V}_{\widehat{\alpha}^{-1}(A^{\widetilde{J}}_{st})},s\geq0)$ converges in distribution as $t\rightarrow0$ to a 1-dimensional symmetric $\alpha$-stable process.
  \item The family of processes $(t^{-1/2}\theta^{V}_{\widehat{\alpha}^{-1}(\exp(st))},s\geq0)$ converges in distribution as $t\rightarrow\infty$ to a 1-dimensional Brownian motion multiplied by $\sqrt{r(\alpha)}$, with
      $$r(\alpha)=\frac{\alpha \ 2^{-1-\alpha/2}}{\pi} \int_{\mathbb{C}} |z|^{-2-\alpha} |\phi(1+z)|^{2} dz,$$
      where $dz$ stands for the Lebesgue measure on $\mathbb{C}$ and for every complex number $z\neq0$, $\phi(z)$ denotes
      the determination of its argument valued in $\left(\right.-\pi,\pi\left.\right]$.
\end{enumerate}
\end{theo}
\vspace{30pt}
{\noindent\textbf{Acknowledgements}} \\
The author would like to thank the anonymous referees for many useful comments and suggestions that improved the quality of this paper.
He is indebted to Professor R.A. Doney for the invitation
at the University of Manchester as a Post Doc fellow where he prepared most of this work and for useful remarks and
to Professor P. Greenwood for motivating this work, for many interesting discussions and for pointing out reference \cite{Gar83}.
The author is also grateful to Professor M. Yor (to the memory of whom this work is dedicated)
for funding from his IUF grant the author's postdoctoral stay at the University of Manchester
and for numerous comments throughout the preparation of this article.



\end{document}